\documentclass[12pt]{amsproc}

\usepackage{amsmath,amssymb,amscd}
\usepackage{layout}
\usepackage{graphicx}

 \theoremstyle{plain}
\newtheorem{theorem}{Theorem}

\newtheorem{corollary}{Corollary}

\def\Z{\mathbb Z}

\def\Q{\mathbb Q}

\def\C{\mathbb C}
\def\F{\mathbb F}
\def\K{\mathbb K}

\def\det{\operatorname{det}}

\def\Gal{\operatorname{Gal}}

\def\mod{\operatorname{mod}}

\def\sgn{\operatorname{sgn}}

\def\SL{\operatorname{SL}}

\begin{document}

\title{Quadratic reciprocity in a finite group}
\author[W. Duke and K. Spears]{William Duke \and Kimberly Spears}\thanks{The research of the second author
was supported by UC LEADS}


\maketitle



\section{Introduction}\label{intro}

The law of quadratic reciprocity is a gem from number theory. In
this article we show that it has a natural generalization to an
arbitrary finite group. Our treatment relies on concepts and
results known at least  100 years ago. See \cite[Chapter I]{Cu}
for a beautiful exposition of the nineteenth century algebra and
number theory we will take as known.

The multiplicative group $(\Z/p\Z)^*$ of reduced residue classes
modulo an odd prime $p$ is a cyclic group of (even) order $p-1$.
Thus it has a unique character of order 2. This character has a
natural pull-back to $\Z$, which is a real Dirichlet character
$(\mod p)$, called the Legendre symbol $\big(\frac{\cdot}{p}
\big)$. By convention, $\big(\frac{a}{p} \big)= 0$ if $p\, |\, a$
and otherwise we have that $\big(\frac{a}{p} \big)=1$ if and only
if $a$ is a square modulo $p$.

In 1872 Zolotarev \cite{Zo} gave an interpretation of the Legendre
symbol $\big(\frac{a}{p} \big)$ that is less well-known; it gives
the sign of the permutation of the elements of $G= \Z/p\,\Z$
induced by multiplication by $a$, provided $p\nmid a$. To see
this, first observe that this recipe defines a character on
$(\Z/p\,\Z)^*$. Furthermore, if it is not trivial, this character
must have order 2 and hence be the Legendre symbol. But it is not
trivial since a generator of $(\Z/p\,\Z)^*$ induces a
$(p-1)$-cycle, which is an odd permutation.
 Motivated by this
observation, we will define in (\ref{def}) below a quadratic
symbol for any finite group $G$.

The classical law of quadratic reciprocity  states that for $q
\neq p$ another odd prime,
\begin{equation}\label{q1}
\big( \tfrac{q}{p} \big)= (-1)^{\frac{p-1}{2}\frac{q-1}{2}}\big(
\tfrac{p}{q} \big), \;\;\;
     \big( \tfrac{-1}{p} \big)=
(-1)^{\frac{p-1}{2}} \;\;\; \textrm{and}\;\;\; \big( \tfrac{2}{p}
\big)= (-1)^{\frac{p^2-1}{8}}.
\end{equation}
 This was first proven by Gauss
in 1796 when he was nineteen years old.  By 1818 he had published
six proofs. The ideas behind his sixth proof \cite{Ga} (see
\cite[p.19]{Cu}), based on the Gauss sum, led to proofs of
quadratic reciprocity using the arithmetic of cyclotomic fields
and the Frobenius automorphism, which was introduced in 1896
\cite{Fr1}. We will combine this classical technique with another
invention of Frobenius  from 1896 \cite{Fr2}, the character table,
to prove a law of reciprocity for the quadratic symbol for any
finite group $G$. A corollary of our result, given in \S \ref{qs},
implies classical quadratic reciprocity when $G=\Z/p\,\Z$ and also
extends Zolotarev's observation to any group of odd order.

\section{The Kronecker symbol}

Before explaining this generalization, let us restate the law of
quadratic reciprocity in one formula by introducing the Jacobi and
Kronecker symbols. The Jacobi symbol simply extends the Legendre
symbol to $\big(\frac{\cdot}{n} \big)$ for arbitrary odd $n \in
\Z^+$ by multiplicativity, so that if $n=p_1 \cdots p_r$ is the
factorization of $n$ into (not necessarily distinct) primes,
$$
\big(\tfrac{a}{n} \big)= \prod_{k = 1}^r \big(\tfrac{a}{p_k}
\big).
$$
 A discriminant is a
non-zero\footnote{We include the possibility that $d$ is a square,
which is usually disallowed.}  integer $d$ with $d \equiv 0
\;\textrm{or} \;1 \;(\mod 4)$.  For a discriminant $d$, the
Kronecker symbol $\big(\frac{d}{\cdot} \big)$ further extends the
Jacobi symbol via
\begin{equation}\label{kro}
\left( \frac{d}{2} \right)=
\left\{%
\begin{array}{ll}
    0, & \hbox{if $d$ is even;} \\
    1, & \hbox{if $d \equiv 1 (\mod 8)$;} \\
    -1, & \hbox{if $d \equiv 5 (\mod 8)$} \\
\end{array}%
\right.
\end{equation}
and by letting $\big(\frac{d}{-1} \big)$  be the sign of $d$. The
value of $\big(\frac{d}{a} \big)$ is then defined for all integers
$a$ by multiplicativity, where we set $\big(\frac{d}{0} \big)=0$
unless $d=1$, in which case $\big(\frac{1}{0} \big)=1$. By means
of these extensions, the law of quadratic reciprocity (\ref{q1})
takes an elegant form for $n$ positive and odd and any integer
$a$:
\begin{equation}\label{q2}
    \left( \frac{a}{n} \right)= \Big( \frac{n^*}{a} \Big),
\end{equation}
where $n^*=(-1)^\frac{n-1}{2}\,n$.   Note that $n^*$ is a
discriminant if $n$ is odd.

\section{The quadratic symbol for a finite group}\label{qs}

Let $G$ be a finite group of order $n$. An integer $a$ that is
prime to $n$ acts as a permutation on the $m$ conjugacy classes
$C_1=\{1\}, C_2, \dots ,C_m$ of $G$ by sending each element $g$ to
$g^a$. Define the quadratic symbol for $G$ at any integer $a$ by
\begin{equation}\label{def}
\left( \frac{a}{G} \right)=\left\{%
\begin{array}{ll}
    0, & \hbox{if $(a, n) \neq 1$;} \\
    1, & \hbox{if this permutation is even;} \\
    -1, & \hbox{if this permutation is odd.} \\
\end{array}%
\right.
\end{equation}
It is easy to see that $\left( \frac{\cdot}{G} \right)$ defines a
real Dirichlet character $(\mod n)$.\footnote{In fact, it is
defined modulo the lcm of the orders of all elements of $G$.}
Zolotarev's observation from the Introduction is that the
quadratic symbol for $G = \Z/p\,\Z$ with an odd prime $p$ is the
Legendre symbol:
\begin{equation}\label{zolo}
\left( \frac{a}{G} \right)=\Big(\frac{a}{|G|} \Big).
\end{equation}

For any group $G$,  a conjugacy class $C$ is said to be real if
$C^{-1}=C$ and complex otherwise. Here $C^{-1}$ denotes the image
of $C$ under $g \mapsto g^{-1}$. Clearly the complex conjugacy
classes occur in pairs $C \neq C^{-1}$ with $|C|=|C^{-1}|.$ Let us
order the conjugacy classes so that the first $r_1$ are real. Thus
$m=r_1+2r_2$ where $r_2$ is half the number of complex conjugacy
classes. Define
\begin{equation}\label{disc}
d=d(G)= (-1)^{r_2}|G|^{r_1} \prod_{j=1}^{r_1} |C_j|^{-1}.
\end{equation}
This is a nonzero integer since for any conjugacy class $C$ we
have that $|G|/|C|=|C_G(g)|$, where $C_G(g)$ is the centralizer of
any $g \in C$ \cite[p.42]{Cu}. It is clear that $d$ is divisible
by $n=|C_G(1)|$ and has the same prime divisors as $n$. We call
$d$ the discriminant of $G$, a name that is justified by the first
statement of our main result.

\begin{theorem}\label{thm}
Let $G$ be a finite group $G$ with discriminant $d$ as defined by
(\ref{disc}). Then $d \equiv 0\textrm{ or} \; 1 \;(\mod 4)$, and
for any integer $a$
\begin{equation}\label{main}
    \left( \frac{a}{G} \right)= \Big( \frac{d}{a} \Big).
\end{equation}
In particular, $\left( \frac{\cdot}{G} \right)$ is trivial if and
    only if $d$ is a square.
\end{theorem}
 In case $G$ has odd order we have the following direct
generalization of classical quadratic reciprocity (\ref{q2}):
\begin{corollary}\label{cor}
If $G$ has odd order $n$ then $d=n^*$ and, for any integer $a$,
\begin{equation}\label{odd}
    \left( \frac{a}{G} \right)= \Big( \frac{n^*}{a} \Big).
\end{equation}
Also, $\left( \frac{\cdot}{G} \right)$ is  trivial if and only if
$n$ is a square.
\end{corollary}
It follows from (\ref{odd}) and (\ref{q2}) that Zolotarev's result
(\ref{zolo}) holds for any group $G$ of odd order.

\section{Proofs}

 The character table of $G$ is
the $m \times m$ matrix  \cite[p.59]{Cu}:
\begin{equation}\label{chart}
M= \begin{pmatrix}\chi_1(C_1) & \dots & \chi_1(C_m)\\
\vdots&\ddots &\vdots\\
\chi_m(C_1) & \dots & \chi_m(C_m)
\end{pmatrix}.
\end{equation}
where $\chi_1=1, \chi_2 \dots, \chi_m$ are the irreducible (over
$\C$) characters of $G$ and we use the convention that
$\chi(C)=\chi(g)$ for any $g \in C$. By the (second) orthogonality
relations we have
\begin{equation}\label{sec}
^t\bar{M}M = \begin{pmatrix}|G||C_1|^{-1} & \dots & 0 \\
\vdots &\ddots &\vdots\\
0 & \dots & |G||C_m|^{-1}
\end{pmatrix},
\end{equation}
a diagonal matrix.  Since $\chi(C^{-1})= \bar{\chi}(C)$ for any
character $\chi$ and any conjugacy class $C$, it is easy to see
that
\begin{equation}\label{sgn}
    \det \bar{M} = (-1)^{r_2} \det M
\end{equation}
and hence by (\ref{disc}) arrive at the identity
\begin{equation}\label{detm}
(\det M)^2 = \ell^2 d
\end{equation}
for some $\ell \in \Z^+.$

Each entry $\chi_i(C_j)$ of $M$ is an algebraic integer in the
cyclotomic field $\Q(\zeta_n)$, where $\zeta_n= e^{2\pi i/n}$. Now
$\Q(\zeta_n)$ is a Galois extension of $\Q$ whose Galois group is
isomorphic to $(\Z/n\,\Z)^*$ by the map $ \sigma_a \mapsto a$,
with $\sigma_a \in \Gal(\Q(\zeta_n)/\Q)$ acting on $\zeta_n$ by
$$\sigma_a(\zeta_n) = \zeta_n^a$$ \cite[Theorem 1 p.92]{Sa}.
 Using this, it is not difficult to check that
\begin{equation}\label{inter}
 \sigma_a(\chi(g)) =
\chi(g^a)
\end{equation}
 for any character $\chi$ and any $g \in G$.

To prove the first statement of Theorem \ref{thm}, we apply an
argument used by Schur \cite{Sc2} to prove Stickelberger's theorem
about the discriminant of a number field. Observe that by the
definition of the determinant
$$\det M= \sum \sgn(\rho)\chi_1(C_{\rho(1)})\chi_2(C_{\rho(2)})\dots
\chi_h(C_{\rho(n)}),$$ where the sum is over all permutations
$\rho$ of the integers $\{1,\dots ,n\}$ and where
$\sgn(\rho)=\pm1$ according to whether $\rho$ is even or odd.
 Write this as $A-B$, where $A$ is the sum of the even permutations
and $B$ is the sum of the odd permutations. By (\ref{inter}) both
of the algebraic integers $A+B$ and $AB$ are invariant under the
Galois group and hence are ordinary integers. In particular, by
(\ref{detm})
$$\ell^2d=(A-B)^2=(A+B)^2-4AB \equiv (A+B)^2 \equiv 0,1 \;(\mod 4),
$$
proving the first statement.

 It is apparent from (\ref{chart}) and (\ref{inter}) that
\begin{equation}\label{frobac}
\sigma_a(\det M)= \left( \frac{a}{G} \right) \det M
\end{equation}
and so by (\ref{detm}) we have
\begin{equation}\label{sqrtd}
\sigma_a(\sqrt{d})= \left( \frac{a}{G} \right)\sqrt{d}.
\end{equation}
Since $\left( \frac{\cdot}{G} \right)$ is a character modulo $n$,
to prove (\ref{main}) it is enough to show it for $a=p$ with $p
\nmid n$ and for $a=-1.$ If $p \nmid n$ we use the Frobenius
automorphism $\sigma_p$. It has the fundamental property that $p$
splits completely in any subfield of $\Q(\zeta_n)$ if and only if
$\sigma_p$ fixes that subfield point-wise \cite[p.91]{Sa}. Thus
$p$ splits in $\K=\Q(\sqrt{d})$ if and only if $\sigma_p
(\sqrt{d})=\sqrt{d}.$

We can always write the discriminant $d$ as an integer square
multiple of the discriminant $d_\K$ of $\K$, called a fundamental
discriminant.  Furthermore, $p$ splits in $\K$  if and only if
$\big(\frac{d_\K}{p} \big)=\big(\frac{d}{p} \big)=1$ \cite[p.
77]{Sa}. Thus we have from (\ref{sqrtd}) that
$$\left( \frac{p}{G}
\right)=\Big( \frac{d}{p} \Big).$$  It is obvious from (\ref{sgn})
and (\ref{disc}) that
\begin{equation}\label{sgn2}
\Big( \frac{-1}{G} \Big)=(-1)^{r_2}=\Big( \frac{d}{-1} \Big),
\end{equation}
finishing the proof of (\ref{main}). That $\big( \tfrac{\cdot}{G}
\big)$ is nontrivial if $d$ is not a square follows easily from
Dirichlet's theorem on primes in progressions. Thus Theorem
\ref{thm} is proven.

Suppose now that $G$ has odd order $n$. Burnside \cite[\S 222
p.294]{Bu} observed that $C_1$ is the only real conjugacy class.
To see this, suppose that $g$ is in a real conjugacy class, so $
h^{-1}gh= g^{-1}$ for some $h$. Then $h^{-2}gh^2= g$ so $h^2 \in
C_G(g)$. Since $n$ is odd, the order of $h$ is odd, say $2\ell+1$.
It follows that $h=(h^2)^{\ell+1}$, so that $h \in C_G(g).$ Thus
$g = g^{-1}$ so $g=1$, since $g$ has odd order.

Since $r_1=1$, it is clear from (\ref{disc}) that $  d
=(-1)^\frac{m-1}{2} n.$ By the first statement of Theorem
\ref{thm} we must have
$$d=(-1)^{\frac{n-1}{2}}n=n^*,$$
since $n$ is odd.\footnote{A stronger result discovered by
Burnside \cite[p.295]{Bu} is that $n \equiv m \; (\mod 16).$} The
last statement of Corollary \ref{cor} follows from that of Theorem
\ref{thm} since if $n$ is odd then $n^*$ is a square if and only
$n$ is a square.

\section{Some examples}

Let us compute the discriminants of some groups with even order.
 Suppose first that $G$ is abelian and that the subgroup of $G$ consisting of 1 and the elements of
 order 2 has order $2^{\,t}$.  Then  $r_1=2^{\,t}$, so
$$
d= (-1)^{\frac{n-2^t}{2}}\;n^{2^t}.
$$
It follows that for $G$ abelian of even order $n$, $\left(
\frac{\cdot}{G} \right)$ is nontrivial if and only if $4\,|\,n$
and $t=1$, in which case for $(a,n)=1$ we have
$$\left( \tfrac{a}{G}
\right)=(-1)^{\frac{a-1}{2}}.$$ The condition $t=1$ holds if $G$
is cyclic, for instance.

In general, if $G$ has only rational characters then it follows
easily from (\ref{inter}) that $\big( \tfrac{\cdot}{G} \big)$ is
the trivial character and hence that $d$ is a square. This holds
in particular for the symmetric group $G=S_k$, where one can also
explicitly compute $d.$

On the other hand, it is not difficult to produce nonabelian
groups with only real characters and with a nontrivial quadratic
symbol. Consider, for example, the family of simple groups given
by $G= \SL(2,\F_{q})$ for $q=2^r$ with $r>1$. By \cite[p.134 (=
p.247)]{Sc1} we have that $n= q(q^2-1)$, that $m=r_1=q+1$ and that
$$d= q^2(q+1)(q^2-1)^{q/2},$$ which is a square if and only if
$r=3$. If $r=2$ we have that $G=A_5$ and $d_\K =5.$ For $r=16$,
$d_\K = 2^{16}+1=65537$, a prime.

We would like to thank Jeffrey Stopple for his comments.


\end{document}